\documentclass[11pt]{amsart}

\usepackage{amssymb}

\usepackage{graphicx}

\newtheorem{theorem}{Theorem}

\theoremstyle{definition}

\theoremstyle{remark}

\newtheorem{example}[theorem]{Example}

\newcommand{\TryPackage}[3]{\IfFileExists{#1.sty}{\usepackage{#1}#2}{#3}}
\TryPackage{mathrsfs}{\renewcommand{\mathcal}{\mathscr}}{%
     \TryPackage{eucal}{}{}}

\newcommand{\lto}{\longrightarrow}

\newcommand{\la}{\lambda}

\newcommand{\Ga}{\Gamma}

\newcommand{\La}{\Lambda}

\newcommand{\ZZ}{{\mathbb Z}}

\newcommand{\CC}{{\mathbb C}}

\newcommand{\SLC}{{SL(2, {\mathbb C})}}

\newcommand{\tr}{\operatorname{\it tr}}

\newcommand{\sm}{{\smallsetminus}}

\begin{document}
\title[Non-triviality of the $M$-degree of the $A$-polynomial]{Non-triviality of the $\boldsymbol{M}$-degree \\of the $\boldsymbol{A}$-polynomial}
 

\author{Hans U. Boden}
\address{Mathematics \& Statistics, McMaster University, Hamilton, Ontario, L8S 4K1 Canada}
\curraddr{}
\email{boden@mcmaster.ca}
\thanks{The author was supported by a grant from the Natural Sciences and Engineering Research Council of Canada.}

\subjclass[2010]{Primary: 57M27, Secondary: 57M25, 57M05}
\keywords{knot; A-polynomial; character variety.}


\begin{abstract}
This note gives a proof that the $A$-polynomial of any nontrivial knot in $S^3$ has nontrivial $M$-degree.\end{abstract}

\maketitle

Given a nontrivial knot $K$ in $S^3$, in this note we prove that its $A$-polynomial satisfies $\deg_M A_K(M,L) \neq 0$. 
This result is established by an argument similar to the one used by Boyer--Zhang and Dunfield--Garoufalidis in their papers on  nontriviality of the $A$-polynomial \cite{BZ3, DG}. 
The question about $\deg_M A_K(M,L)$  comes up naturally in the study of the $SL(2,\CC)$ invariant  of knots, and the examples in \cite{M, IMS} show that the $A$-polynomial of a knot can have a nontrivial $L$-factor.   
The result of \cite{BZ3,DG}  implies that the $A$-polynomial detects the unknot, but it is extremely difficult to compute. Apart from the low-crossing knots, whose $A$-polynomials can be found at \cite{ChLi} thanks to M. Culler and J. Hoste, there are only a handful of infinite families of knots for which the $A$-polynomial is known, namely the torus knots \cite{CCGLS}, the twist knots and several sub-families of generalized twist knots \cite{HS, P}, and a sub-family of  pretzel knots \cite{GM, TY}.  
Presumably, the $M$-degree of $A_K(M,L)$ is easier to compute than the full $A$-polynomial,  and the main result in this paper  shows that it is sufficient for detecting  unknottedness.

We begin by introducing the notation for $\SLC$ character varieties and the $A$-polynomial.
For a finitely generated group $\Ga$, we use $R(\Gamma)$ to denote the
space of representations $\varrho\colon \Ga \lto \SLC$ and note that $R(\Ga)$ is a complex affine
algebraic set (see \cite{CS} or \cite{Sh}). The {\sl character}  of a representation $\varrho$ is
the function $\chi_\varrho\colon \Ga \lto \CC$ defined by 
$\chi_\varrho(g)=\tr(\varrho(g))$ for $g \in \Ga$. We denote the set
of characters by $X(\Ga)$ and note that it also admits the structure of
 a complex affine algebraic set (see again \cite{CS} or \cite{Sh}).
For a manifold $M$, we write
$R(M)$ for the variety of $\SLC$ representations of $\pi_1 (M)$ and
$X(M)$ for its character variety.

Given a knot $K$, let $N_K = S^3 \sm \tau(K)$ be its complement and choose  a standard meridian-longitude pair $(\mu,\la)$ for $\pi_1(\partial N_K)$. Let $\La \subset R(\partial N_K)$ be the subset of diagonal representations, i.e.
$$\La = \{ \varrho:\pi_1(\partial N_K) \lto \SLC \mid \text{$\varrho(\mu)$ and  $\varrho(\la)$ are both diagonal matrices} \},$$
and define the eigenvalue map  $\La \lto \CC^* \times \CC^*$ by sending $\varrho \mapsto (u,v) \in \CC^* \times \CC^*$, where $$\varrho(\mu) = \begin{pmatrix}u& 0 \\ 0 &u^{-1}\end{pmatrix} \text{ and }\varrho(\la) = \begin{pmatrix}v& 0 \\ 0 &v^{-1}\end{pmatrix}.$$
This map identifies $\La$ with  $\CC^* \times \CC^*$, and the natural projection $p:\La \lto X(\partial N_K)$ is a degree 2, surjective, regular map.

The  inclusion $i:\partial N_K \lto N_K$ induces a map on $\pi_1$ which in turn gives rise to the map $i^*:X(N_K) \lto X(\partial N_K)$ which is regular. We define $V \subset X(\partial N_K)$ to be the Zariski closure of the union of the image $i^*(X_j)$ over each component $X_j \subset X(N_K)$ for which $i^*(X_j)$ is one-dimensional, and we set
$D \subset \CC^2 $ to be the Zariski closure of the algebraic curve $p^{-1}(V) \subset \La,$ where we  identify $\La$ and  $\CC^* \times \CC^*$ via the eigenvalue map.
Thus $D \subset \CC^2$ is a plane curve, and we define $A_K(M,L)$ to be its defining polynomial. 
It is well-defined up to sign by the requirement that its coefficients lie in $\ZZ$ and have greatest common divisor equal to one, and that it have no repeated factors. The component of abelian characters gives rise to the factor $L-1$, and for the unknot $\bigcirc$, one can easily show that
$A_\bigcirc(M,L) = L-1.$ By the main theorem in \cite{BZ3, DG},  if $K$ is a nontrivial knot in $S^3$, then $A_K(M,L) \neq L-1.$ 
 
\begin{theorem} \label{nontrivM}
If $K$ is a nontrivial knot in $S^3$, then its $A$-polynomial satisfies 
$\deg_M A_K(M,L) \neq 0.$
\end{theorem}

\noindent
We point out that,  just like nontriviality of the $A$-polynomial, the above statement is easy to verify  in the case of a torus or hyperbolic knot. 

\begin{proof}
Assume $K$ is a non-trivial knot and set $N_K = S^3 \sm \tau(K)$. 
Suppose to the contrary that $\deg_M A_K(M,L) =0$. Then by Theorem 2.9 of \cite{CL2} (see also Proposition 10.3 of \cite{BZ2}), it follows that $$A_K(M,L) = (L-1) f_1(L) \cdots f_k(L),$$ where 
$f_1,\ldots, f_k$ are distinct cyclotomic polynomials. Choose positive integers $d_1, \ldots, d_k$ so that the roots of $f_i(L)$ have order $d_i$ for each $i=1,\ldots, k$.  

In \cite{KM}, P. Kronheimer and T. Mrowka establish the existence of an infinite family $\{ \varrho_n :\pi_1(M) \to SU(2)\}$ of irreducible representations, where $\varrho_n$ extends over $1/n$ surgery on $K$. Moreover, by Claim 2.1 of \cite{DG}, it follows that the restrictions of $\varrho_n$ to $\pi_1(\partial M)$ are all distinct as points  in $X(\partial M)$.  Let  $(u_n, v_n) \in \CC^* \times \CC^*$ denote the point corresponding to the restriction of $\varrho_n$ to $\pi_1(\partial M)$, and by the surgery condition, it follows that $u_n v_n^n = 1.$ Further, all but finitely many of the pairs $(u_n, v_n)$ satisfy the $A$-polynomial $A_K(M,L)$. Fix $d = d_1 \cdots d_k$ and consider the infinite family $\{\varrho_{nd} \mid n = 1,2,3, \ldots \}.$ Since all but finitely many of the pairs $\{(u_{nd}, v_{nd})\}$ satisfy the $A$-polynomial, it follows that for at least one $n$, say $v_{nd}$, is a root of   $(L-1) f_1(L) \cdots f_k(L)$.
But every root $\xi$ of this polynomial satisfies $\xi^d = 1,$ and this together with the surgery condition implies that $u_{nd} = (v_{nd})^{nd} = 1.$ As $\varrho_{nd}$ is an $SU(2)$ representation, this implies $\varrho_{nd}(\mu) = I,$ and since the meridian is a normal generator for the knot group $\pi_1(M)$, this  implies that $\varrho_{nd}:\pi_1(M) \to \SLC$ is trivial, which contradicts irreducibility of the representation $\varrho_{nd}$.  
 \end{proof}

There is an example due to T. Mattman \cite{M} that helps to put this result into context. In order to describe it, we  recall a definition of a refinement of the $A$-polynomial due to S. Boyer and X. Zhang. For every  one-dimensional component $X_i$ of $X(N_K)$ with image $i^*(X_j)$ also one-dimensional, they define an $A$-polynomial $A_{X_j}(M,L)$ in \cite{BZ2}. Basically, their idea is to apply the above construction to the given component $X_j$ and to incorporate the degree of $i^*|_{X_j}$ as the multiplicity of the resulting factor.
Set
$$\widehat{A}_K(M,L) =  A_{X_1}(M,L) \cdots A_{X_n}(M,L),$$ the product taken over all one-dimensional components $X_j$ of $X(N_K)$ with one-dimensional restriction $i^*(X_j)$. For a small knot $K$ in $S^3$, since by \cite{CS} no component of $X(N_K)$ can have dimension larger than one, it follows that $A_K(M,L)$ and $\widehat{A}_K(M,L) $ coincide up to repeated factors.  

\begin{example}
In \cite{M}, Mattman shows that the character variety of the $(-2,3,-3)$ pretzel knot   contains a $0$-curve, and so consequently its $A$-polynomial has a nontrivial $L$-factor.
(According to \cite{IMS}, this holds for all $(-2,p,-p)$ pretzel knots provided $p$ is odd.)
In other words, there exists a component $X_j$ of $X(N_K)$ with $A$-polynomial satisfying $\deg_M A_{X_j}(M,L)=0$. This shows that Theorem \ref{nontrivM} is not true in the more general setting of the Boyer-Zhang $A$-polynomials. 
\end{example}

We outline an alternative approach to proving Theorem \ref{nontrivM}, which involves reducing to a case where one can use the identity  
\begin{equation} \label{A-form}
A_K(\pm 1,L)=\pm L^a(L-1)^b(L+1)^c \text{ for integers } a,b,c
\end{equation}
(see section 2.8 of \cite{CCGLS} and Corollary 4.5 of \cite{CL1}). This formula holds for most but not all knots;  among the computations, one finds that $9_{29}$ and $9_{38}$ give counterexamples. The Newton polytopes for both knots have a vertical edge, and indeed Corollary 4.5 in \cite{CL1} includes the hypothesis that $K$ is a small knot, meaning $N_K$ does not contain an incompressible surface that is not boundary parallel.

The idea for reducing to this case is due to M. Culler, and it is based on the observation that Equation (\ref{A-form}) holds for any knot whose Newton polytope does not have a vertical edge. The main result of \cite{CL2} implies $A_K(\pm 1,L)$ is monic, and the argument in section 2.8 of \cite{CCGLS} applies to show that Equation (\ref{A-form}) holds (cf. Corollary 4.5 in \cite{CL1}). Moreover, by \cite[Theorem 3.4]{CCGLS}, we see that  Newton polytope has a vertical edge if and only if $\infty$ is a strongly detected boundary slope. In other words, if $\infty$ is not a strongly detected boundary slope, then Equation (\ref{A-form}) holds.

Suppose then that $\infty$ is a strongly detected boundary slope. Then the vertical edge of the Newton polytopes ensures that $\deg_M A_K(L,M) \neq 0,$ and so we are reduced  to the case $\infty$ is not a strongly detected boundary slope. By Equation (\ref{A-form}), if $A_K(L,M)$ were constant in $M$, then \cite{BZ3, DG} implies $A_K(L,M) = (L-1)(L+1).$ The proof is completed by arguing as in Theorem \ref{nontrivM} with $k=1$ and the single cyclotomic polynomial $f_1(L) =L+1$.

\begin{figure}[t]
\begin{center}
\leavevmode\hbox{}
\includegraphics[scale=0.26]{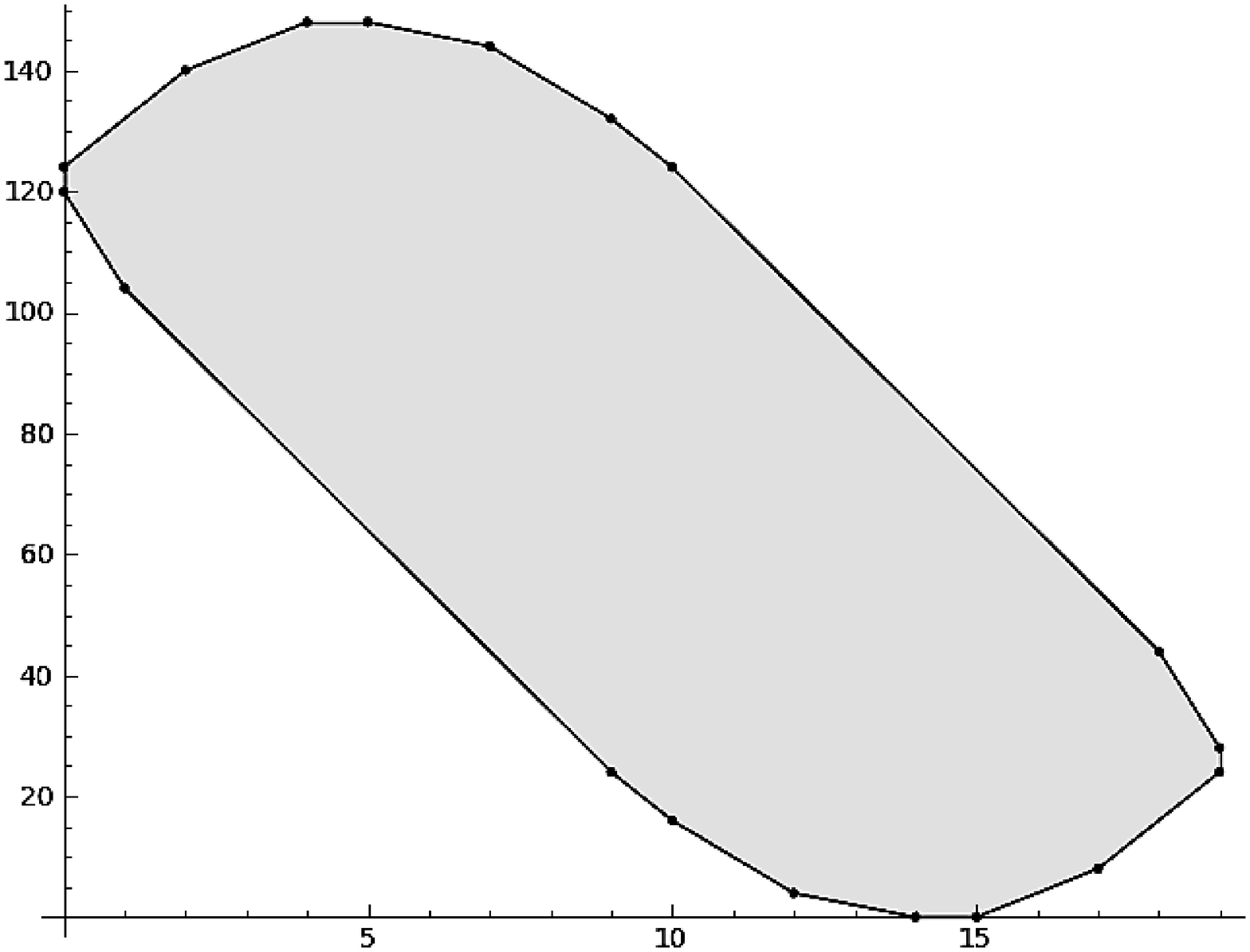} 
\includegraphics[scale=0.26]{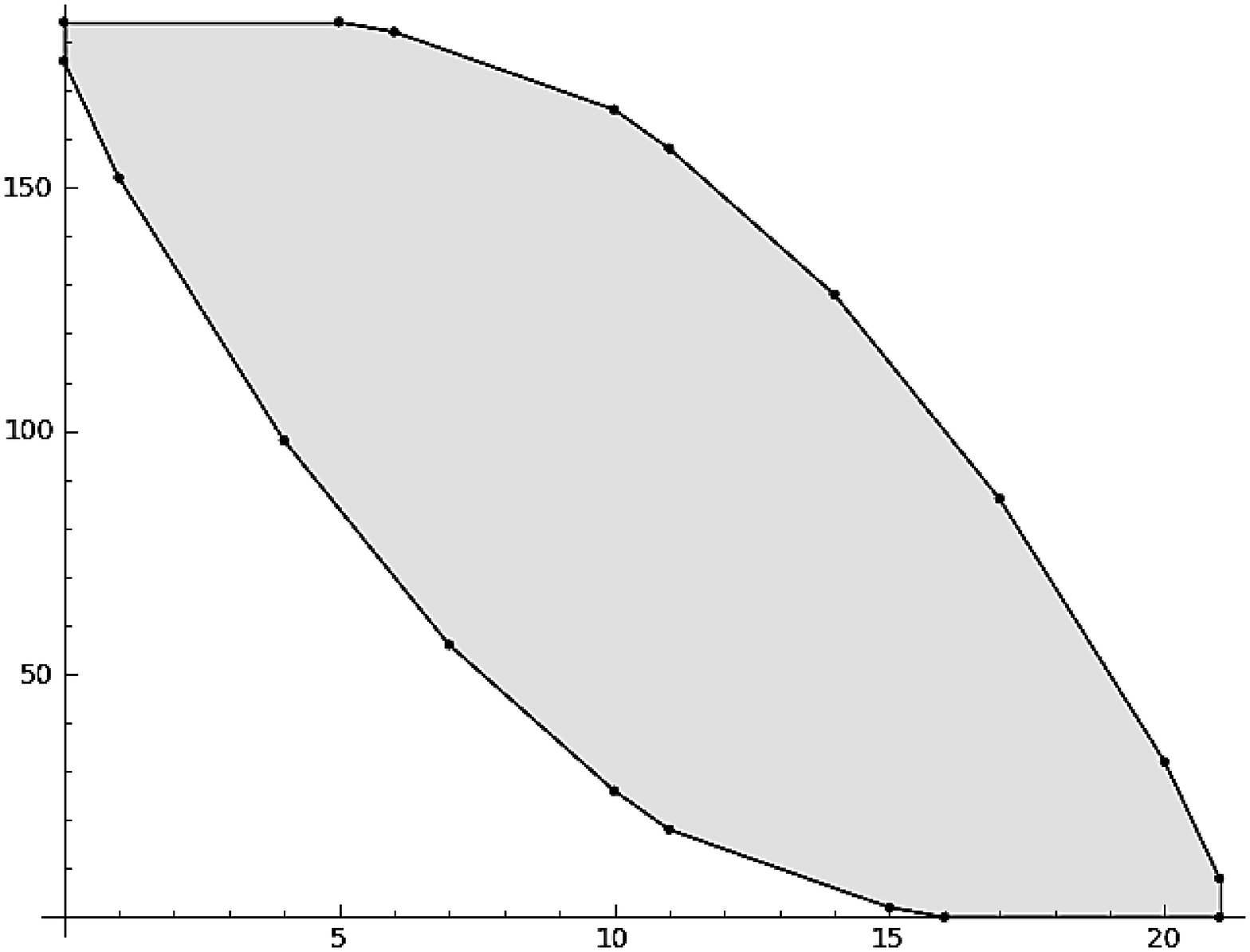} 
\caption{Newton polytopes for $9_{29}$ and $9_{38}$} \label{newt}
\end{center}
\end{figure}

\medskip 
Before concluding, we briefly indicate the relationship between this work and questions about the $\SLC$ Casson knot invariant $\la'(K)$ (see \cite{C, BC} for more details). 
 Using
the surgery theorem for the $\SLC$ Casson invariant, one can show that for any small knot $K$ in $S^3$, $\la'(K) \neq 0 $ if and only if $\deg_M A_K(M,L) \neq 0.$ Essentially, both quantities can be interpreted in terms of Culler-Shalen seminorms $\| \la \|_{CS}$ of the longitude, and in \cite{BC} we
use this to prove nontriviality of $\la'(K)$ for all small knots. 

Whether $\la'(K) \neq 0$ holds more generally is not immediately clear; in fact the invariant is at present only defined for small knots (cf. \cite{C}).  The general question whether  $\la'(K)$ detects the unknot is closely related to the question whether  $\deg_M \widehat{A}_K(M,L) \neq 0$ for all nontrivial knots. Note that this does not follow from Theorem \ref{nontrivM}. In fact, despite the positive results in \cite{BZ3, DG}, it is an open problem to show nontriviality of $\widehat{A}_K(M,L)$ for all nontrivial knots $K$. In terms of the character variety, this problem is related to the question of existence of irreducible $\SLC$ representations $\varrho$ of the knot group whose characters $\chi_\varrho$ lie on a one-dimensional algebraic component
$X_j \subset X(N_K)$ whose image $i^*(X_j)$ in $X(\partial N_K)$ is also one-dimensional.

 \medskip \noindent
{\it Acknowledgements.} The author is extremely grateful to Marc Culler for his numerous contributions to this project. The author would also like to thank Steve Boyer, Cindy Curtis, Stavros Garoufalidis, and Thomas Mattman for helpful discussions. Finally, the author thanks the Max Planck Institute for Mathematics in Bonn for its support. 
 
\bibliographystyle{amsplain}

\end{document}